\documentclass[12pt,final]{article}

\usepackage[body={6in,9 in}, top=1 in, left=1in, bottom=1in, right=1in, nohead]{geometry}
\usepackage{amsfonts, amsmath, amssymb, amsthm}
\usepackage{systeme}
\usepackage{latexsym}
\usepackage{enumerate}
\usepackage{verbatim}
\usepackage{graphicx}
\usepackage[all,2cell,ps]{xy}
\usepackage{setspace} 
\usepackage{titlesec}
\usepackage[titles]{tocloft}
\usepackage{ragged2e}
\usepackage{indentfirst}
\usepackage[utf8]{inputenc}
\usepackage{hyperref}
\usepackage{cleveref}
\usepackage{float}
\usepackage[shortlabels]{enumitem}
\usepackage{ kantlipsum}
\usepackage{listings}
\usepackage{xcolor}
\usepackage{mathtools}

\usepackage{colorprofiles}
\usepackage[a-2b,mathxmp]{pdfx}[2019/02/27]
\hypersetup{pdfstartview=}

\RaggedRight

\graphicspath{ {./images/} }

\newtheoremstyle{newthm} 
     {2pt}
     {2pt}
     {}
     {}
     {\bfseries}
     {.}
     {.5em}
     {}

\theoremstyle{newthm}
\newtheorem{theorem}{Theorem}[section]
\newtheorem{thm}{Theorem}[section]
\newtheorem*{theorem*}{Theorem}
\newtheorem{lem}{\textbf{Lemma}}
\newtheorem{lemma}{\textbf{Lemma}}

\newtheorem{dff}[lem]{Definition}
\newtheorem{defn}[lem]{Definition}

\newtheorem{cor}[lem]{Corollary}

\crefname{Theorem}{Theorem}{Theorems}
\crefname{Proposition}{Proposition}{Propositions}
\crefname{Definition}{Definition}{Definitions}
\crefname{Corollary}{Corollary}{Corollaries}
\crefname{Example}{Example}{Examples}
\crefname{Remark}{Remark}{Remarks}
\crefname{Conjecture}{Conjecture}{Conjectures}
\crefname{Lemma}{Lemma}{Lemmas}
\crefname{prop}{Proposition}{Propositions}
\newcommand{\FF}{\mathbb F}
\newcommand{\ZZ}{\mathbb Z}

\newcommand{\CC}{\mathbf C}

\newcommand{\inv}{^{-1}}

\newcommand{\bet}{\beta}

\newcommand{\alp}{\alpha}
\newcommand{\al}{\alpha}

\newcommand{\eps}{\epsilon}

\newcommand{\varedgeII}{\Xi}
\newcommand{\vartreeNumII}{\Upsilon}
\newcommand{\vartreeNumIII}{\Phi}
\newcommand{\edgeSetVar}{\mathcal{F}}

\newcommand{\mat}[1]{\begin{bmatrix} #1 \end{bmatrix}}

\newcommand{\rar}{\rightarrow}

\newcommand{\nin}{\not\in}
\newcommand{\ovl}[1]{\overline{#1}}

\newcommand{\tm}{\times}

\newcommand{\Tr}{\operatorname {Tr}}

\begin{document}
\pagenumbering{arabic}

 \title{Dot Product Bounds in Galois Rings}
 \author{David Crosby}

 \maketitle

\noindent \textbf{ABSTRACT}

\smallskip

\begin{singlespace}
 We consider a problem related to the Erd\H os Unit Distance Conjecture: How often can a single dot product configuration or a multiple dot product configuration occur over a Galoi Ring? We also find a bound on an inverse vector matrix multiplication problem.
\end{singlespace}

\noindent\textbf{KEYWORDS}: Dot Product, Erd\H os Unit Distance Conjecture, Galois Rings

\section*{Acknowledgements} I want to thank Cameron Wickham for suggesting this problem and Steven Senger for his many helpful conversations.

\setcounter{tocdepth}{2}
\tableofcontents




\section{Introduction}\label{introChap}

In 1946, Paul Erd\H os asked how many unit distances may occur in the plane with $n$ points? Erd\H os proved that $n$ points in $R^2$ determines $\Omega(n^{1/2})$ many unit distances and conjectured that there are $\Omega(n^c)$ for any $c<1$ (see \cite{erdosConj}). The best known result is $\Omega(n/\log(n))$ by Guth and Katz \cite{guthKatz} in 2015. For $R^d$, Erd\H os proved that $n$ points determines $\Theta(n^{2/d})$ unit distances and in 2008 Solymosi and Vu proved the bound of $\Omega(n^{2/d-2/(d(d+2))})$ \cite{solymosi}.
The idea of the unit distance problem has since been extended to various rings by using a distance like functional. We see this in Covert, Iosevich and Pakianathan's 2011 work \cite{1105} where they study integers modulo a prime power and also Iosevich and Rudnev's 2008 work \cite{siminal} over finite fields.

 A common extension to this question is the single dot product problem: ``How often can a specified dot product occur between $n$ points in some ambient space?''. This has been studied in \cite{avgHyper}, \cite{upBnds}, \cite{pairsfieldsAndRings}, and more, where different bounds and configurations are studied.
 

\subsection{Focus of the article}
This article will focus on various configurations of the dot product problem in the ambient space of Galois rings. 
Galois rings are a set of finite rings that contain both all of $F_{p^l}$ and $\ZZ_{p^l}$. They are constructed as
$R_{e,k}\cong\ZZ_{p^e}[x]/(f)$ where $f$ is of degree $k$ and irreducible in $\ZZ_p[x]$ under the mod $p$ map. 
The configurations that we consider are single dot products, pairs of dot products, and general forest configurations. We further show an application of the forest configuration leads to a bound on the number of solutions to the matrix equation $aA=b$ where $b$ is a given column and $a$ and $A$ are our variables.

Galois rings are finding applications in computer networking through coding theory. For more on Galois rings in coding theory see \cite[\S 8.1.3]{biniFlam}.
Another reason to study finite rings is that they provide a test bed for complicated analysis problems as integrals are always convergent. The prototypical example of this is Zeev Dvir's paper ``On the size of Kakeya sets in finite fields'' \cite[Theorem 2]{ZeevDvir}, where he proved the Kakeya Conjecture (relates to subsets of $\FF_{p^l}$ that contain a line in every direction) via simple means.
We study forest configurations as they encompass all dot product configurations that contain no cycles.

\subsection{Results}

Our first result follows in form the proof by Covert, Iosevich, and Pakianathan \cite[Theorem 1.3.2]{1105}. There they study the Erd\H os distance problem over the integers modulo an odd prime to a power. 
We define our counting function $\nu$, which counts how many pairs of points give a specified dot product.
\begin{defn}\label{nuDef}
   Let $E \subset R_{e,k}^d$. Let $p,e,k,d$ be given natural numbers with $p$ prime, $e\geq 5$, $d\geq 2$, and $k\geq 1$.  
   Define $$\nu(t) = \left|\{(x,y) \in E\tm E: x\cdot y = t\}\right|.$$
\end{defn}

Covert et al. in \cite{1105} show that $\nu(t) = |E|^2/p^l+R(t)$, where $$|R(t)|< l|E|p^{l\left( \frac{d-1}{2} (2-1/l)\right)}<|E|^2/{p^l}$$ whenever $|E| > lp^{ld-d/2+1/2}$. 

\begin{figure}[hb]
      \center{
      \includegraphics{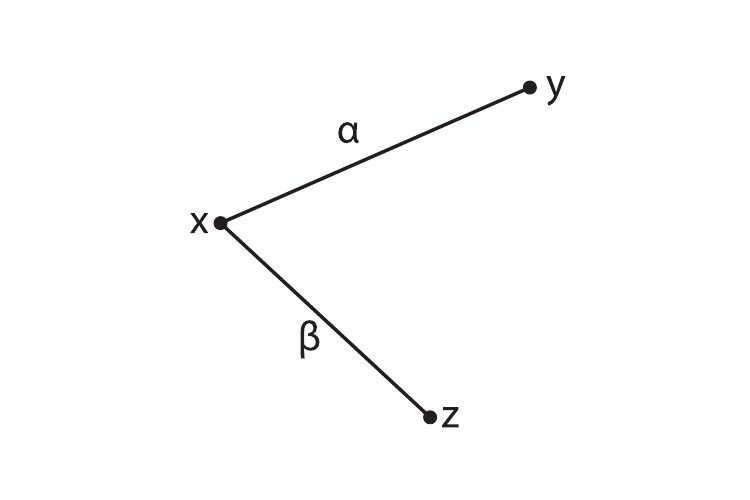}
      }
      \caption{Pair of Dot Products}
      \label{figPairDots}
   \end{figure}

Our first result gives that the number of pairs of points that give a specified dot product must be no more than twice the expected value whenever the size of our point set is large enough. 
\begin{thm}\label{singleDotClean}
   Let $E \subset R_{e,k}^d$. Let $p,e,k,d$ be given natural numbers with $p$ prime, $e\geq 5$, $d\geq 2$, and $k\geq 1$.  
 Let $\nu(t) = |\{(x,y) \in E\tm E: x\cdot y = t\}|$.
   Then 
   $$\nu(t) \leq 2|E|/p^{ek}$$ for any $t\in R_{e,k}$ whenever $|E|\geq \sqrt{6+3e}p^{dek-dk/2+ek/2+k/2}$.
\end{thm}
   As $|E|\leq|R^d|=p^{dek}$, this result is non-trivial when $d\geq e+1+\frac{\log_p(6+3e)}{k}$. 
   Our second results gives that the number of triplets of points $(x,y,z)$ that give specified dot products ($x\cdot y = \al,\, y\cdot z = \bet$) must be no more than twice the expected value whenever the size of our point set is large enough. Graphically, this is given in Figure \ref{figPairDots}.
    We define a new counting function to count these triplets of points.

\begin{dff}\label{PiDef}
   Let $d\geq 3,\ e\geq 5,\ k\geq 1$. Let $E\subset (R_{e,k})^d$ and suppose that $\al,\bet \in R_{e,k}$. We define 
   $$\Pi_{\al,\bet}(E) = \{ (x,y,z) \in E^3 : x\cdot y = \al, x\cdot z = \bet\}.$$
\end{dff}
The proof follows and extends the paper ``Pairs of Dot Products in Finite Fields and Rings''  by David Covert and Steven Senger \cite{pairsfieldsAndRings}, where they consider both $\FF_{p^l}$ and $\ZZ_{p^l}$. 
Covert and Senger find 
$$\left|\Pi_{\al,\bet}(E)\right| = \frac{|E|^3}{p^{2l}}(1+o(1)) $$ whenever $|E|> cp^{l\left( \frac{d+1}{2}\right)}$ for $E\subset \FF_{p^l}^d$ and $|E|>cp^{dl-d/2+1}$ when $E\subset \ZZ_{p^l}^d$ for some constant $c$.

\begin{thm}\label{pairsDotClean}
   Let $E \subset R_{e,k}^d$. Let $d\geq 3,\ e\geq 5,\ k\geq 1$. Let $E\subset (R_{e,k})^d$ and suppose that $\al,\bet \in R_{e,k}$. We have the bound $$|\Pi_{\al,\bet}(E)|\leq \frac{2|E^3|}{p^{2ek}}$$ whenever 
   $|E| \geq \sqrt{35e/8}p^{dek-dk/2+ek/2+k/2}.$
\end{thm}

\begin{figure}[ht]
   \center{
\includegraphics{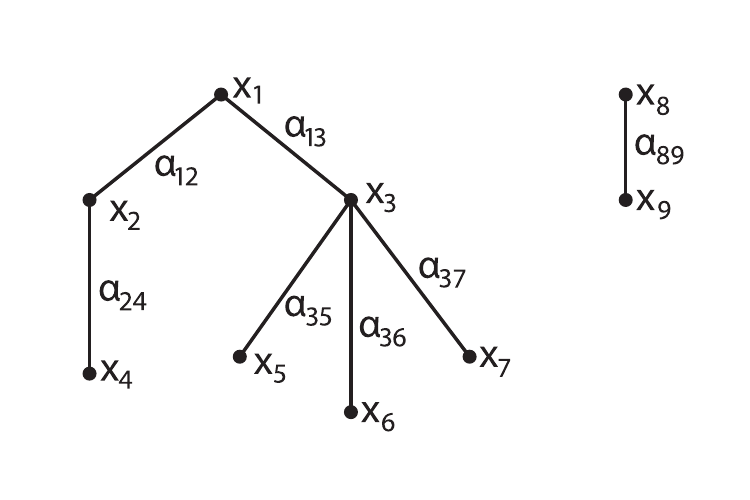}
}
\caption{Forest of elements and dot products}
\label{figForest}
\end{figure}

This result is non-trivial when $d\geq e+1+\log_p(35e/8)/k$.
Our third result is similar, but for tree configurations. 
We extend the results of \cite{kChainsFq} by Blevins, Lynch, Senger, and the author from chains over $\ZZ_p^l$ and $\FF_p^l$ to forests over Galois rings.
We create a third counting function as follows: 
\begin{dff}\label{PiTreeDef}
   Let $E \subset R_{e,k}^d$. Let $T$ be a forest with vertices $V$ with $|V|=m$ where each vertex represents an element of $E$ and edges $\edgeSetVar$ with $|\edgeSetVar|=n$  where each edge $(e_i,e_j)$ represents that we have specified $e_i\cdot e_j = \al_{ij}$ where $\al = \{\al_{\eps}\}_{\eps\in \edgeSetVar}$ is given. An example graph is visualized in Figure \ref{figForest}. We define
   $\Pi_\al(T) = \{ x\in E^m : \forall \al_{ij}\in \al, x_i\cdot x_j = \al_{ij}\}.$
\end{dff}
\begin{theorem}\label{trees}
   Let $d\geq 3,\ e\geq 5,\ k\geq 1,\ n\geq 3$. Let $E \subset R_{e,k}^d$. 
   Let $T$ be a forest with vertices $V$ with $|V|=m$ where each vertex represents an element of $E$ and edges $\edgeSetVar$ with $|\edgeSetVar|=n$  where each edge $(e_i,e_j)$ represents that we have specified $e_i\cdot e_j = \al_{ij}$ where $\al = \{\al_{\eps}\}_{\eps\in \edgeSetVar}$ is given. Then $$\Pi_\al(T) \leq \frac{2|E|^m}{p^{nek}} $$ whenever 
   $|E| \geq \max\left\{2n\sqrt{6+3e}p^{dek-dk/2+3ek/2+k/2}, \sqrt{14}ep^{dek-dk/2+nek/2+n/2-ek/2+k/2}\right\}.$
\end{theorem}
   When $n\geq 4$, the bound on $|E|$ reduces to $|E| \geq \sqrt{14}enp^{dek-dk/2+nek/2+n/2-ek/2+3k/2}.$ This is nontrivial (with $n\geq 4$) when $d  \geq (n-1)e+3+\frac{n+\log_p(14en^2)}{k}$.
An application of Theorem \ref{trees} gives a bound on how many pairs of row-vectors and matrices there are that multiply to a specified column-vector.
\begin{theorem}\label{invRowMatProb}
   Let $d\geq 3,\ e\geq 5,\ k\geq 1,\ n\geq 4$. Let $E \subset R_{e,k}^d$. 
   Let $aA=b$ be an equation where $a\in E$, the columns of $A$ are in $E$, and $b$ is a given row-vector of length $n$.  Whenever $|E|> en\sqrt{14}p^{dek-dk/2+nek/2+n/2-ek/2+k}$,  the number of solutions $$S=\left|\left\{(a,A)\mid a\in E, A\in E^{n},\  aA=b\right\}\right|$$ has 
$ S \leq \frac{2|E|^{n+1}}{p^{nek}}. $
\end{theorem}
We now compare this bound to the known number for when $e=k=1$ and $E=R^d_{e,k}=\ZZ_p^d$. From Hodges 1966 work \cite{HodgesModPa}, where he gives an expression for the number of represntations of a bilinear form ($\bmod\ p^a$), we are able to find that the number of such solutions $S$  is $p^{dn-n}(p^{d}-1)$ which is less than twice the expected value of $\frac{2p^{dn+d}}{p^{n}}$.
The proof for Theorem \ref{invRowMatProb} works by showing a similarity of structures. So 
should there be a bound on the matrix-matrix multiplication problem $AB=C$, then its bound should give a bound on a corresponding dot product configuration with cycles.

\subsection{Background}

For proofs of the below statements see \cite{biniFlam}, \cite{holdman}, and \cite{mcDonald}.

\begin{dff}[$\chi_{e,k}$]\label{canAddChar}
   
   The canonical additive character $\chi_{e,k}:R_{e,k}\rar \CC^\tm$ is defined as $$\chi_{e,k}(z)=e^{2\pi i\Tr_{e,k}(z)/p^e}.$$
\end{dff}

\begin{lem}[Orthogonality]\label{orthog}
   Let $e,k$ be given. For any $a\in R$,
   $$\sum_{z\in R} \chi(az)=\begin{cases} |R|, & a=0 \\ 0, & a\neq 0\end{cases}$$
\end{lem}

\begin{thm}\label{tDef}
   Let $\bet$ be a generator of the subgroup of $R_{e,k}^\tm$ isomorphic to $\ZZ_{p^k-1}$. Let $T_{e,k}= \{0,1,\bet,\ldots, \bet^{p^k-2}\}$.  Every $z\in R$ has a unique $p$-adic representation $$z=z_0+pz_1+\cdots + p^{e-1}z_{e-1},\quad z_i\in T_{e,k}.$$
\end{thm} 

\begin{dff}
We will use the notation $[p^i]$ to mean $(p^i)\setminus (p^{i+1})$. 
\end{dff}

\begin{lem}\label{pBrakStruc}
   Any $s\in [p^i]$ has the form $s=p^iu$ where $u$ is the uniquely determined unit of the form $u_1+pu_2+\cdots+p^{e-i-1}u_{e-i-1}$ where $u_i\in T_{e,k}$. 
\end{lem}

\begin{cor}\label{multNonPos}
   If $a\in R_{e,k}$ is given, then
   $$ \sum_{z\in p^iR_{e,k}}\chi_{e,k}(az) = \sum_{w\in R_{e-i,k}}\chi_{e-i,k}(\rho_i(a)w).$$
   Where $\rho_i(a_0+pa_1+\cdots+ p^{e-1}a_{e-1})=a_0+pa_1+\cdots+p^{e-1-i}a_{e-1-i}$ with $a_0+pa_1+\cdots+ p^{e-1}a_{e-1}$ being the $p$-adic expansion of $a$.
\end{cor}
\begin{proof}
    By definition of $\chi_{e,k}$,  $\Tr$, and characteristic of a ring, this is immediate. Let $(az)_0+\cdots+p^{e-1}(az)_{e-1}$ be the $p$-adic expansion of $az$. Note,
   \begin{align*}
      \sum_{z\in p^i R_{e,k}}\chi_{e,k}(az) 
   \end{align*}
   \begin{align*}
                                           &= \sum_{z\in R_{e,k}}\exp\left(2\pi i\left[p^iaz+\left(p^i(az)_0^p+\cdots+ p^{e-1+i}(az)^p_{e-1}\right)+\cdots +\left(p^i(az)_0^{p^k}+\ldots p^{e-1+i}(az)^{p^k}_{e-1}\right)\right]/p^e\right) \\
                                           &= \sum_{z\in R_{e,k}}\exp\left(2\pi i\left[az+\left(p^i(az)_0^p+\cdots+p^{e-1-i}(az)^p_{e-1}\right)+\cdots +\left((az)_0^{p^k}+\cdots+ p^{e-1}(az)^{p^k}_{e-1}\right)\right]/p^{e-i}\right) \\
                                           &= \sum_{w\in R_{e-i,k}}\chi_{e-i,k}(\rho_i(a)w).
   \end{align*}
\end{proof}

\begin{lem}\label{unitsLeq0} 
    Let $e,k$ be given and $R=R_{e,k}$. If $n$ is a given natural number and $0<n<e$, then
   $$\sum_{z\in R^\tm} \chi(p^nz) \leq 0.$$
\end{lem}
\begin{proof}
   Note that 
   \begin{align*}
      \sum_{z\in R^\tm} \chi(p^nz) &= \sum_{z\in R}\chi(p^nz) - \sum_{z\in (p)}\chi(p^nz) \\
                                   &= I+II. \\
   \end{align*}
   By orthogonality (Lemma \ref{orthog}), $I=0$.
   So $\sum_{z\in R^\tm}\chi(p^nz)=-\sum_{z\in (p)}\chi(p^nz)$. Thus, when $n+1=e$, $\sum_{z\in R^\tm}\chi(p^nz)=-|(p)|=-|pR_{e,k}|$. 
   In the case that $n+1\neq e$, we apply Corollary \ref{multNonPos},
   \begin{align*}
      \sum_{z\in R^\tm} \chi(p^nz) & = -\sum_{pz\in pR_{e,k}}\chi(p^{n+1}z) \\
                                    & = -\sum_{w\in R_{e-1,k}}\chi_{e-1,k}(p^nw). \\
   \end{align*}
   By orthogonality in $z$, this is $0$ (as $n\neq e-1$). Thus, $\sum_{z\in R^\tm} \chi(p^nz) \leq 0$ when $n<e$. Further, $$\sum_{z\in R^\tm}\chi(p^nz) = \begin{cases} -|pR_{e,k}| &,\ n+1=e \\ 0 &,\ n+1 < e \end{cases}$$
   
\end{proof}

For more properties of character sums over Galois rings, see F. Shuqin and H. Wenbao \cite{charSums}.

\section{Proof of Single Dot Product Result}\label{secSingDotProd}



\noindent Theorem \ref{singleDotClean} is an application of the more technical Theorem \ref{singleDot} whose proof we delay until later.
\begin{thm}\label{singleDot}
   Let $e,k$ be given natural numbers greater than 0. Let $E\subset R_{e,k}^d$. Let $\nu(t) = |\{(x,y) \in E\tm E: x\cdot y = t\}|$. For any $t\in R$, 
   $$\nu(t)< |E|^2/p^{ek} + D(t).$$ 
   Further, the discrepancy $D(t)$, has $D(t)<|E|^2/p^{ek}$ whenever 

      $$|E|  \geq p^{ek}\sum_{i=0}^{e-1} \sqrt{p^{dek+dik-ik}\left(2p^{-ik-k}+1+p^{-ek}\right)}.$$
\end{thm}
   We now prove Theorem \ref{singleDotClean}.
\begin{proof}
      Label 
   \begin{equation}\label{singDotCleanFsetup}
      F=p^{ek}\sum_{i=0}^{e-1} \sqrt{p^{dek+dik-ik}\left(2p^{-ik-k}+1+p^{-ek}\right)}.
   \end{equation}
   By Theorem \ref{singleDot}, we know that 
   \begin{align*}
      |E| & \geq F. \\
      \intertext{Examining $F$, we see }
      F  & = p^{ek}\sum_{i=0}^{e-1} \sqrt{p^{dek+dik-ik}\left(2p^{-ik-k}+1+p^{-ek}\right)} \\
         &=p^{ek+dek/2}\sum_{i=0}^{e-1} \sqrt{p^{dik-ik}\left(2p^{-ik-k}+1+p^{-ek}\right)}. \\
   \end{align*}
   So, 
   \begin{align*}
      \frac{F}{p^{ek}p^{dek/2}} &\leq \sum_{i=0}^{e-1} p^{dik/2-ik/2}\sqrt{2p^{-ik-k}+1+p^{-ek}}. \\
         \intertext{For getting a cleaner bound for Theorem \ref{singleDot}, we now apply Cauchy--Schwarz, which \allowbreak loosens our bound to}
      \frac{F}{p^{ek}p^{dek/2}} & \leq \sqrt{\sum_{i=0}^{e-1}(p^{dik/2-ik/2})^2 \sum_{i=0}^{e-1}\sqrt{2p^{-ik-k}+1+p^{-ek}}^2} \\
                                & \leq \sqrt{\sum_{i=0}^{e-1}p^{dik-ik} \sum_{i=0}^{e-1}2p^{-ik-k}+1+p^{-ek}}. \\
      \intertext{This is a sum of geometric series or constants. Thus,}
      \frac{F}{p^{ek}p^{dek/2}} & \leq \sqrt{\frac{1-p^{dek-ek}}{1-p^{dk-k}}\left( 2p^{-k}\frac{1-p^{-ek}}{1-p^{-k}}   +e+ep^{-ek}\right)}. \\
                                  \intertext{The above fraction in parenthesis is maximized when $e$ grows large, $k=1$, and $p=2$, giving $ \frac{1-p^{-ek}}{1-p^{-k}} \leq 2$. 
                               Ergo, }
      \frac{F}{p^{ek}p^{dek/2}} & \leq \sqrt{\frac{1-p^{dek-ek}}{1-p^{dk-k}} \left(4p^{-k}+e+ep^{-ek}\right)} \\
                                & \leq \sqrt{\left(\frac{1}{1-p^{dk-k}} - \frac{p^{dek-ek}}{1-p^{dk-k}}\right) \left(4p^{-k}+e+ep^{-ek}\right)}. \\
                                \intertext{Let $C=1/(1-p^{dk-k})$. We may bound $-\frac{p^{dek-ek}}{1-p^{dk-k}}$ by $2p^{dek-ek-dk+k}$ as $\frac{-1}{1-B}\leq \frac{2}{B}$ for any $B\geq 2$. This gives,}
      \frac{F}{p^{ek}p^{dek/2}} &\leq \sqrt{(C + 2p^{dek-ek-dk+k})\left(4p^{-k}+e+ep^{-ek}\right)}. \\
                                \intertext{Since $e\geq 5$, $k\geq 1$, and $p\geq 2$, we have $4p^{-k} +e +ep^{-ek}\leq 69/32+e $. So,}
      \frac{F}{p^{ek}p^{dek/2}} &\leq \sqrt{(C + 2p^{dek-ek-dk+k})\left(69/32+e\right)}. \\
      \intertext{As $C=1/(1-p^{dk-k})$, $p^{dk-k}\geq 2$, $|C|\leq 1$, and $p^{dek-ek-dk+k}\geq 2$, we have $C+2p^{dek-ek-dk+k} \leq (|C|/2+2)p^{dek-ek-dk+k}$. So,}
      \frac{F}{p^{ek}p^{dek/2}} &\leq \sqrt{(1/2 + 2)p^{dek-ek-dk+k}\left(69/32+e\right)} \\
                                & \leq \sqrt{6+3e}p^{dek/2-dk/2-ek/2+k/2}. \\
   \end{align*}

   Thus 
   \begin{equation}\label{singDotCleanFfound}
   F\leq \sqrt{6+3e}p^{dek-dk/2+ek/2+k/2}
   \end{equation}
   and so $\nu(t) \leq 2|E|/p^{ek}$ whenever $$|E|> \sqrt{6+3e}p^{dek-dk/2+ek/2+k/2}.$$
\end{proof}

   As $|E|\leq|R^d|=p^{dek}$, this result is non-trivial when
   \begin{align*}
      p^{dek} & \geq \sqrt{6+3e}p^{ dek-dk/2+ek/2+k/2}, \\
      \intertext{which simplifies to}
      p^{dk/2} & \geq \sqrt{6+3e}p^{ek/2+k/2} \\
      dk & \geq ek+k+\log_p(6+3e). \\
   \end{align*}
   So whenever $d\geq e+1+\frac{\log_p(6+3e)}{k}$, then we may take a proper subset of $R^d$ for Theorem \ref{singleDotClean}.
   We now get to the heart of this paper; the technical result for the number of pairs of points in a Galois ring that have a specified dot product. The proof of Theorem \ref{singleDot} is given below.
\begin{proof}
   Recall that by Lemma \ref{orthog}, $\sum_{r\in R}\chi(ra)=0$ if $a\neq 0$ and $p^{ek}$ if $a=0$. Since we are after $x\cdot y = t$, we examine $\sum_{s\in R}\sum_{x,y\in E} \chi(s(x\cdot y -t))$. This sum gives out a $p^{ek}$ precisely when $x\cdot y = t$ and 0 when $x\cdot y \neq t$. Thus, by multiplying this sum by $p^{-ek}$ we get the number of $x,y \in E$ such that $x\cdot y = t$.  We rewrite $\nu(t)$ (Definition \ref{nuDef}) as below and split $\nu(t)$ into the following parts 
   \begin{align*}
      \nu(t) &= p^{-ek}\sum_{s\in R}\sum_{x,y\in E}\chi(s(x\cdot y - t)) \\
             &= \nu_0(t)+\cdots +\nu_e(t),\\
   \end{align*}
   where $$\nu_i(t) = p^{-ek}\sum_{s\in [p^i]}\sum_{x,y\in E}\chi(s(x\cdot y))\chi(-st).$$

   For $\nu_e(t)$, we have $$\nu_e(t) = p^{-ek}\sum_{x,y\in E} \chi(0(x\cdot y))\chi(-0t) = |E|^2/p^{ek}.$$ 
   The discrepancy is 
   \begin{equation}\label{singDdefEq}
      D(t) = \sum_{i=0}^{e-1} \nu_i(t) =p^{-ek}\sum_{s\in R\setminus \{0\}}\sum_{x,y\in E}\chi(s(x\cdot y - t)) .
   \end{equation}

   We now examine $\nu_i(t)$ for some $i \neq e$. Recall,
   
   \begin{align*}
      \nu_i(t) &= p^{-ek}\sum_{s\in [p^i]}\sum_{x,y\in E}\chi(s(x\cdot y))\chi(-st).\\
      \intertext{We will  prepare to use Cauchy-Schwarz by recognizing that }
      \nu_i(t) &= p^{-ek}\left(\sum_{x\in E}\left(1\right)\sum_{y\in E}\sum_{s\in [p^i]}\chi(s(x\cdot y))\chi(-st)\right) \\
               &= \left(\sum_{x\in E}\left(1\right)\sum_{y\in E}p^{-ek}\sum_{s\in [p^i]}\chi(s(x\cdot y))\chi(-st)\right). \\
               \intertext{ By applying Cauchy-Schwarz with $a_x = 1$ and $b_x = \sum_{y\in E}p^{-ek}\sum_{s\in [p^i]}\chi(s(x\cdot y))\chi(-st)$,} 
      |\nu_i(t)|^2 & \leq p^{-2ek}\left(\sum_{x\in E} 1\bar{1}\right)\left(\sum_{x\in E}\sum_{y\in E}\sum_{s\in [p^i]}\chi(s(x\cdot y))\chi(-st)\ovl{\sum_{y'\in E}\sum_{s'\in [p^i]}\chi(s'(x\cdot y'))\chi(-s't)}\right).\\
      \intertext{Because $z\bar z\geq 0$ for any $z\in \CC$, we may dominate the second sum of $x\in E$ by $x\in R^d$, obtaining,}
      |\nu_i(t)|^2& \leq p^{-2ek}|E|\sum_{x\in R^d} \sum_{y,y' \in E}\sum_{s,s'\in [p^i]}\chi(s(x\cdot y-t))\chi(-s'(x\cdot y'-t)). \\
      \intertext{  
         By Lemma \ref{pBrakStruc}, $s\in [p^i]$ is the same as $p^iu\in [p^i]$ where $u\in T_{e,k}$.
 Thus as $s,s'\in [p^i]$, we  let $s=p^iu$ and $s'=p^iv$. Using $\ovl0$ to be the zero element of $R_{e,k}^d$ this gives, }
   \end{align*}
   \begin{equation}\label{nui^2Eq}
                   |\nu_i(t)|^2 \leq p^{-2ek}|E|\sum_{x\in R^d} \sum_{y,y'}\sum_{p^iu,p^iv\in [p^i]}\chi(p^i(uy-vy')\cdot x)\chi(p^it(u-v)) 
    \end{equation}
    \begin{align*}
                   & = p^{-2ek}|E| \sum_{y\in E}\left(\sum_{\substack{y'\in E,\ p^iu,p^iv\in [p^i] \\ p^i(uy-vy')=\ovl{0}}}\sum_{x\in R^d}\chi(p^i(uy-vy')\cdot x)\chi(p^it(u-v)) \right. \\
                   & \qquad+ \left. \sum_{\substack{y'\in E,\ p^iu,p^iv\in [p^i] \\ p^i(uy-vy')\neq \ovl{0}}}\sum_{x\in R^d}\chi(p^i(uy-vy')\cdot x)\chi(p^it(u-v))\right). \\
                   \intertext{By orthogonality in $x$, whenever $p^i(uy-vy')=\bar{0}$ we get a factor of $p^{ek}$ for each component of $x$.  When $p^i(uy-vy')\neq 0$, we get a factor of zero. Thus, }
                   |\nu_i(t)|^2 & \leq |E|p^{dek-2ek}\sum_{\substack{y,y'\in E \\ p^i(uy-vy')=\ovl{0} \\ p^iu,p^iv\in [p^i]}} \chi(p^it(v-u)). \\
                   \intertext{We now split this sum into}
                   |\nu_i(t)|^2&\leq I+II,
   \end{align*}
   where $I$ has $u=v$ and $II$ has $u\neq v$.

   \begin{lem}\label{singILem}
      For $I$, we have
      $|I| \leq |E|^2p^{dek+ikd-ek-ik}.$
   \end{lem}
   \begin{proof}
   For $I$, we have
   \begin{align*}
      I & =|E|p^{dek-2ek}\sum_{\substack{y,y'\in E \\ p^i(uy-vy')=\ovl{0} \\ p^iu=p^iv\in [p^i]}} \chi(p^it(v-u)), \\
      \intertext{which as $u=v$,}
      I & =|E|p^{dek-2ek}\sum_{\substack{y,y'\in E \\ p^iu(y-y')=\ovl{0} \\ p^iu\in [p^i]}} \chi(p^it(0)). \\
       \intertext{ Let $E(y)=1$ when $y\in E$ and $E(y)=0$ when $y\nin E$. As $u$ is a unit, $p^iu(y-y')=\ovl 0$ is the same as $p^i(y-y')=\ovl 0$. As $\chi(p^it0)=1$, we have}
      I & =|E|p^{dek-2ek}\left|[p^i]\right|\sum_{\substack{y,y'\in R^d \\ p^i(y-y')=\ovl{0} }} E(y)E(y') \\
       & =|E|p^{dek-2ek}\left(p^{(e-i)k}-p^{(e-i-1)k}\right)\sum_{\substack{y,y'\in R^d \\ p^i(y-y')=\ovl{0} }} E(y)E(y'). \\
   \end{align*}

   Note $$\sum_{\substack{y,y'\in R^d \\ p^i(y-y')=\ovl0}}E(y)E(y')$$ is the count of $(y,y')\in E^2$ such that $y-y' \in (p^{e-i})^d$. This we may bound above by taking an arbitrary $y\in E$ and seeing that there are, at most, $|p^{e-i}R|=p^{ik}$  choices for each component of $y'$. This gives $$\sum_{\substack{y,y'\in R^d \\ p^i(y-y')=\ovl0}}E(y)E(y') \leq \sum_{y\in E}\sum_{\substack{y'\in R^d \\p^i(y-y')=\bar 0}} 1 \leq |E|p^{ikd}.$$
   Thus, $$|I| \leq|E|p^{dek-2ek}\left(p^{(e-i)k}-p^{(e-i-1)k}\right)|E|p^{ikd} \leq  |E|^2p^{dek+ikd-ek-ik}.$$
   \end{proof}

   For $II$ we have $v\neq u$. So,
   \begin{align*}
      II & = |E|p^{dek-2ek}\sum_{y\in E}\sum_{\substack{y'\in E,\ p^iu,p^iv\in [p^i] \\ p^i(uy-vy')=\ovl0 \\ u\neq v}}\chi(p^it(v-u)). \\
      \intertext{Let $v=b$ and $a=u/v$. This gives us,}
      II & = |E|p^{dek-2ek}\sum_{y\in E}\sum_{\substack{y'\in E,\ p^iu,p^iv\in [p^i] \\ p^ib(ay-y')=\ovl0 \\ u\neq v}}\chi(p^itb(1-a)) \\
         &= II_u + II_n, \\
   \end{align*}
   where $II_u$ has $1-a$ being a unit and $II_n$ has $1-a$ being a non-unit.

   \begin{lem}\label{singIIuLem}
      For $II_u$ we have,
      $$|II_u| \leq |E|^2p^{dek+dik-2ik-k} + |E|^2p^{dek+dik-ik}.$$    
   \end{lem}
   \begin{proof}
   For $II_u$, we have that as $1-a$ is a unit,
   \begin{align*}
      II_u &= |E|p^{dek-2ek}\sum_{y\in E} \sum_{\substack{y'\in E,\ p^iu,p^iv\in [p^i] \\ p^ib(ay-y')=\ovl0 \\ u\neq v \\ 1-a \in R_{e-i,k}^\tm }}\chi(p^itb(1-a)).\\
      \intertext{As $b=v$ and $p^iv\in [p^i]$, $b$ has the form $b=u_1+pu_2+\cdots + p^{e-i-1}u_{e-i-1}$ where $u_i \in T_{e,k}$ by Lemma \ref{pBrakStruc}. So, summing over $p^ib\in [p^i]$ is the same as summing over $b\in R_{e-i,k}^\tm$. 
   As $a=u/v$ and $1-a\neq 0$,  it must be that $u\neq v$ (allowing us to drop it from the restriction on the summation). This gives, }
          |II_u| &\leq |E|p^{dek-2ek}\sum_{y\in E} \sum_{\substack{y'\in E,\ b\in R^\tm_{e-i,k} \\p^ib(ay-y')=\ovl0 \\  1-a \in R_{e-i,k}^\tm }}\chi_{e-i,k}(tb(1-a)).\\
   \end{align*}
   Since $p^ib(ay-y')=\bar 0$, and $b$ is a unit, we must have $ay-y'\in (p^{e-i})^d$, so
   \begin{align*}
          |II_u|  &\leq |E|p^{dek-2ek}\sum_{y\in E} \sum_{\substack{y'\in E,\ b\in R^\tm_{e-i,k} \\ay-y' \in (p^{e-i})^d \\  1-a \in R_{e-i,k}^\tm }}\chi_{e-i,k}(tb(1-a)).\\
   \end{align*}
   As $a$ is a unit, and $ay-y'\in (p^{e-i})^d$, it must be that for each choice of $y$,  that $y'$ is in the same coset of $(p^{e-i})^d$. 
   Being the case that $\chi_{e-i,k}(tb(1-a))$ does not depend on $y$ nor $y'$, 
    we can bound $ay-y'\in (p^{e-i})^d$  by summing $y\in E$ and summing $y'\in -ay+(p^{e-i})^d$. This gives  $|(p^{e-i})^d|=p^{idk}$  choices for $y'$. Pulling out this factor of $p^{idk}$, we get
   \begin{align*}
           |II_u| &\leq |E|p^{dek-2ek+idk}\sum_{y\in E} \sum_{\substack{1-a \in R_{e-i,k}^\tm }}\sum_{b\in R^\tm_{e-i,k}} \chi_{e-i,k}(tb(1-a)).\\
   \end{align*}

   As $b$ sums over $R_{e-i,k}^\tm$, by Lemma \ref{unitsLeq0},  we have that $II_u=0$ when $t\nin (p^{e-i-1})$.
   When $t\in [p^{e-i-1}]$, we have, by Lemma \ref{unitsLeq0},
   \begin{align*}
      |II_u|&\leq \left||E|p^{dek-2ek+idk}\sum_{y\in E} \sum_{\substack{1-a \in R_{e-i,k}^\tm }}-\left|pR_{e-i,k}\right|\right|,\\
      \intertext{which then simplifies to}
      |II_u| & \leq \left|E\right|p^{dek-2ek+idk}\left|E\right|\left|R^\tm_{e-i,k}\right|\left|pR_{e-i,k}\right| \\ 
          & \leq |E|^2p^{dek+dik-2ik-k}. \\
   \end{align*}

   In the final case that $t\in (p^{e-i})$, we have $\chi(tb(1-a))=1$, and so
   \begin{align*}
      |II_u| &\leq |E|p^{dek-2ek+idk}\sum_{y\in E} \sum_{\substack{1-a \in R_{e-i,k}^\tm }}|R_{e-i,k}^\tm|.\\
      \intertext{Which then simplifies as follows:}
      |II_u| & \leq \left|E\right|p^{dek-2ek+idk}\left|E\right|\left|R^\tm_{e-i,k}\right||R_{e-i,k}^\tm| \\ 
           & \leq |E|^2p^{dek+dik-ik}.
   \end{align*}
   Thus, $$|II_u| \leq |E|^2p^{dek+dik-2ik-k} + |E|^2p^{dek+dik-ik}.$$

\end{proof}

\begin{lem}\label{singIInLem}
   For $II_n$, we have
   $$|II_n|  \leq |E|^2p^{dek+dik-2ik-k}.$$
\end{lem}
\begin{proof}
   Recall $v=b$ and $a=u/v$. For $II_n$, we have $1-a$ being a non-unit, hence $p^i(1-a)\in (p^{i+1})$. This requirement makes $u\neq v$ satisfied. Thus,
   \begin{align*}
      II_n & = |E|p^{dek-2ek}\sum_{y\in E}\sum_{\substack{y'\in E,\ p^iv,p^iu\in [p^i] \\ p^ib(ay-y')=\ovl0 \\  p^i(1-a) \in (p^{i+1})}}\chi(p^itb(1-a)). \\
      \intertext{ As summing over $p^iu \in [p^i]$ is the same as summing over $p^iu/v \in [p^i]$, we have, }
      II_n & = |E|p^{dek-2ek}\sum_{y\in E}\sum_{\substack{y'\in E,\ p^ib,p^ia\in [p^i] \\ p^ib(ay-y')=\ovl0 \\  p^i(1-a) \in (p^{i+1})}}\chi(p^itb(1-a)). \\
      \intertext{Recall $p^ia \in [p^i]$ means that $a=u_0+pz_1+\cdots +p^{e-i-1}z_{e-i-1}$ for some unit $u_0$ and $z_j\in T_{e,k}$ by Lemma \ref{pBrakStruc}. As $1-a$ is not a unit, $a=1+pr$ for some $r\in R_{e,k}$. This means that the restriction $p^i(1-a)\in (p^{i+1})$ restricts $p^ia \in [p^i]$ to the subset $p^ia \in p^i+(p^{i+1})$. Thus,}
      II_n & = |E|p^{dek-2ek}\sum_{y\in E}\sum_{\substack{y'\in E,\ p^ib\in [p^i],\ p^ia\in p^i+(p^{i+1}) \\ p^ib(ay-y')=\ovl0  }}\chi(p^itb(1-a)). \\
   \end{align*}
   By letting $c=p^i(1-a)$, we have
   \begin{align*}
           II_n & = |E|p^{dek-2ek}\sum_{y\in E}\sum_{\substack{y'\in E,\ p^ib\in [p^i],\ c\in (p^{i+1}) \\ p^ib(ay-y')=\ovl0 }}\chi(tbc). \\
           \intertext{As $p^ib(ay-y')=\ovl0$ means $y' \in ay+ (p^{e-i})^d$ and no where else do we depend on $y'$, we bound this sum by summing $y'\in (p^{e-i})^d$. This gives,}
           |II_n| & \leq  |E|p^{dek-2ek}\sum_{y\in E}\sum_{\substack{y'\in (p^{e-i})^d,\ p^ib\in [p^i],\ c\in (p^{i+1})  }}\chi(tbc). \\
           \intertext{Which, by Corollary \ref{multNonPos} (supressing $\rho_i$) and pulling out $y'\in (p^{e-i})^d$, is}
           |II_n| & \leq  |E|p^{dek-2ek}\left|(p^{e-i})^d\right|\sum_{y\in E}\sum_{\substack{p^ib\in [p^i]}}\sum_{c\in R_{e-i-1,k}}\chi(tbc). \\
   \end{align*}
   By orthogonality in $c$, Lemma \ref{orthog}, we have $II_n=0$ when $t\nin (p^{e-i-1})$. Otherwise,
   \begin{align*}
           |II_n| & \leq  |E|p^{dek-2ek}\left|(p^{e-i})^d\right|\sum_{y\in E}\sum_{\substack{p^ib\in [p^i]}}\sum_{c\in R_{e-i-1,k}}1. \\
           \intertext{We simplify this to}
           |II_n| & \leq |E|^2p^{dek+dik-2ik-k}.
   \end{align*}
\end{proof}

   Recall,
   \begin{align*}
      |II|& \leq |II_u|+|II_n|, \\
   \intertext{which by Lemmas \ref{singIIuLem} and \ref{singIInLem} give,} 
      |II| & \leq \left(|E|^2p^{dek+dik-2ik-k} + |E|^2p^{dek+dik-ik}\right) + |E|^2p^{dek+dik-2ik-k} \\
          & \leq 2|E|^2p^{dek+dik-2ik-k} + |E|^2p^{dek+dik-ik}. \\
   \end{align*}
   As $|\nu_i(t)|^2 \leq |I|+|II|$ and by our bounds on $II$ and $I$ (Lemma \ref{singILem}), 
    \begin{equation}\label{singDotNuIEq}
       |\nu_i(t)|^2 \leq 2|E|^2p^{dek+dik-2ik-k} + |E|^2p^{dek+dik-ik} + |E|^2p^{dek+ikd-ek-ik}. 
    \end{equation}

   Thus the discrepancy, $D(t)=\left|\sum_{i=0}^{e-1} \nu_i(t)\right|$, has 
   \begin{align*}
      D(t) & \leq \sum_{i=0}^{e-1} \sqrt{2|E|^2p^{dek+dik-2ik-k} + |E|^2p^{dek+dik-ik} + |E|^2p^{dek+ikd-ek-ik}} \\
            & \leq \sum_{i=0}^{e-1} \sqrt{|E|^2p^{dek+dik-ik}\left(2p^{-ik-k}+1+p^{-ek}\right)}. \\
   \end{align*}
   Thus,
   \begin{equation}\label{singDiscrIneq}
      D(t) \leq |E|\sum_{i=0}^{e-1} \sqrt{p^{dek+dik-ik}\left(2p^{-ik-k}+1+p^{-ek}\right)}.
   \end{equation}
   So as we want $\nu_e(t)\geq D(t)$, we must have
   \begin{align*}
      |E|^2/p^{ek} & \geq |E|\sum_{i=0}^{e-1} \sqrt{p^{dek+dik-ik}\left(2p^{-ik-k}+1+p^{-ek}\right)}\\
      |E| & \geq p^{ek}\sum_{i=0}^{e-1} \sqrt{p^{dek+dik-ik}\left(2p^{-ik-k}+1+p^{-ek}\right)}.\\
   \end{align*}
   This concludes the proof of Theorem \ref{singleDot}.
\end{proof}

   We have also proved the following two corollaries.
\begin{cor}\label{dotProdRemNu}
   Let $p$ be prime, $e\geq 5,\,k\geq 1,\, d\geq 2$ be given. With $R=R_{e,k}$, 
   \begin{align*}
        & p^{-2ek}|E| \sum_{y,y'\in E}\sum_{p^iu,p^iv\in [p^i]}\sum_{x\in R^d}\chi(p^i(uy-vy')\cdot x)\chi(p^it(u-v)) \\ 
   & \leq 2|E|^2p^{dek+dik-2ik-k} + |E|^2p^{dek+dik-ik} + |E|^2p^{dek+ikd-ek-ik} 
   \end{align*}
   whenever 
   $$|E|\geq p^{ek}\sum_{i=0}^{e-1} \sqrt{p^{dek+dik-ik}\left(2p^{-ik-k}+1+p^{-ek}\right)}$$
\end{cor}
\begin{proof}
   By Equation \ref{nui^2Eq} we have our left hand side is greater than $|\nu_i(t)|^2$ and by Equation \ref{singDotNuIEq} we have our right hand side is greater than $|\nu_i(t)|^2$. As we only relaxed our bound between the two equations, this follows.
\end{proof}
\begin{cor}\label{dotProdRem}
   Let $p$ be prime, $e\geq 5,\,k\geq 1,\, d\geq 2$ be given. With $R=R_{e,k}$, 
   \begin{align*}
    & p^{-ek}\sum_{s\in R\setminus \{0\}}\sum_{x,y \in E}\chi(s(x\cdot y - t))\\
     & \leq |E|\sum_{i=0}^{e-1} \sqrt{p^{dek+dik-ik}\left(2p^{-ik-k}+1+p^{-ek}\right)} \\
   \end{align*}
   whenever 
   $$|E|\geq p^{ek}\sum_{i=0}^{e-1} \sqrt{p^{dek+dik-ik}\left(2p^{-ik-k}+1+p^{-ek}\right)}$$
\end{cor}
\begin{proof}
   By the definition of $D(t)$ (Equation \ref{singDdefEq}) and Equation \ref{singDiscrIneq} this is immediate.
\end{proof}
\begin{cor}\label{dotProdRemSimp}
   Let $p$ be prime, $e\geq 5,\,k\geq 1,\, d\geq 2$ be given. With $R=R_{e,k}$, 
$$p^{-ek}\sum_{s\in R\setminus \{0\}}\sum_{x,y \in E}\chi(s(x\cdot y - t)) \leq \sqrt{6+3e}|E|p^{dek-dk/2+ek/2+k/2}$$
   whenever 
   $$|E|\geq p^{ek}\sum_{i=0}^{e-1} \sqrt{p^{dek+dik-ik}\left(2p^{-ik-k}+1+p^{-ek}\right)}$$
\end{cor}
\begin{proof}
	By Corollary \ref{dotProdRem}, and the relaxations made between Equations \ref{singDotCleanFsetup} and \ref{singDotCleanFfound} from the proof of Theorem \ref{singleDotClean} this is clear.
\end{proof}


\section{Proofs of Mulitiple Dot Produt Results}
\subsection{Proof of Pairs of Dot Products Result}\label{pairsChap}

\noindent Theorem \ref{pairsDotClean} relies on the more technical Theorem \ref{doubDotProd} whose statement is below and whose proof is given later in this section.
\begin{thm}\label{doubDotProd}
    Let $e,k\geq 1$. Let $d\geq 3$, $E\subset (R_{e,k})^d$ and suppose that $\al,\bet \in R_{e,k}$.  We have the bound $$|\Pi_{\al,\bet}(E)|=\frac{|E^3|}{p^{2ek}}(1+o(1))$$ whenever 
      $$|E| \geq \max\left(p^{ek}\sum_{i=0}^{e-1} \sqrt{p^{dek+dik-ik}\left(2p^{-ik-k}+1+p^{-ek}\right)}, \sqrt{2G}\right)$$
      where 
   \begin{align*}
      G &=p^{dek+2ek}e\left(2p^{-k}\frac{p^{dek-2ek}-1}{p^{dk-2k}-1} + (1+  p^{-ek})\frac{p^{dek-ek}-1}{p^{dk-k}-1}\right).
   \end{align*}

\end{thm}

We now begin the proof of Theorem \ref{pairsDotClean}.
\begin{proof}
   Let $e\geq 5, d\geq 3$, and $k\geq 1$. By 
   Theorem \ref{doubDotProd} 
   we know that  $$|\Pi_{\al,\bet}(E)|=\frac{|E^3|}{p^{2ek}}(1+o(1))$$ is true whenever 
      $$|E| \geq \max\left(p^{ek}\sum_{i=0}^{e-1} \sqrt{p^{dek+dik-ik}\left(2p^{-ik-k}+1+p^{-ek}\right)}, \sqrt{2G}\right)$$ where
   \begin{align*}
      G &= p^{dek+2ek}e\left(2p^{-k}\frac{p^{dek-2ek}-1}{p^{dk-2k}-1} + (1+  p^{-ek})\frac{p^{dek-ek}-1}{p^{dk-k}-1}\right).
   \end{align*}
   By 
   Equations \ref{singDotCleanFsetup} and \ref{singDotCleanFfound}, we know $$ p^{ek}\sum_{i=0}^{e-1} \sqrt{p^{dek+dik-ik}\left(2p^{-ik-k}+1+p^{-ek}\right)}\leq \sqrt{6+3e}p^{dek-dk/2+ek/2+k/2}$$ whenever $e\geq 5,\ d\geq 2$, and $k\geq 1$, which is given. We focus now on getting a nice upper bound on $G$. Recall
   \begin{align*}
      G &=p^{dek+2ek}e\left(2p^{-k}\frac{p^{dek-2ek}-1}{p^{dk-2k}-1} + (1+  p^{-ek})\frac{p^{dek-ek}-1}{p^{dk-k}-1}\right). \\
      \intertext{As $1/(B-1)\leq 2/B$ for any $B\geq 2$, }
      G &\leq 2p^{dek+2ek}e\left(2p^{-k}\frac{p^{dek-2ek}-1}{p^{dk-2k}} + (1+  p^{-ek})\frac{p^{dek-ek}-1}{p^{dk-k}}\right)\\
          &\leq 2p^{2dek+ek-dk+k}e\left(3p^{-ek} + 1\right).\\
       \intertext{Recall, $e\geq 5, p\geq 2, k\geq 1$, so the part in parentheses is maximized when $e=5,\ p=2,$ and $k=1$. Thus, } 
      G &\leq 2p^{2dek+ek-dk+k}e\left(3(2^{-5})  + 1\right) \\
       &\leq p^{2dek+ek-dk+k}e(35/16). \\
   \end{align*}
   Thus, this lemma is true whenever $$|E|\geq \max\left(\sqrt{6+3e}p^{dek-dk/2+ek/2+k/2}, \sqrt{2(35/16)ep^{2dek+ek-dk+k}}\right).$$
   As $\sqrt{35e/8}p^{dek-dk/2+ek/2+k/2} \geq \sqrt{6+3e}p^{dek+ek/2-dk/2+k/2}$ since $e\geq 5$, we have that 
   $$|E| \geq \sqrt{35e/8}p^{dek-dk/2+ek/2+k/2}.$$
\end{proof}

This is non-trivial when $|E|< p^{dek}$, giving 
   $p^{dek} \geq \sqrt{35e/8}p^{dek-dk/2+ek/2+k/2}.$
Hence, this result is non-trivial when $d\geq e+1+\log_p(35e/8)/k$.
   We now begin the proof of Theorem \ref{doubDotProd}. 
\begin{proof}
   For the sake of brevity, we will use $R$ to stand for $R_{e,k}$. Let $\chi$ denote the canonical additive character of $R$ (Definition \ref{canAddChar}).
   Recall that,
   \begin{align*}
      |\Pi_{\al,\bet}(E)|&=|\{(x,y,z)\in E\tm E\tm E : x\cdot y = \alp, x\cdot z = \bet\}|.\\
      \intertext{As we want $x\cdot y = \al$ and $x\cdot z = \bet$, the character sum becomes}
      |\Pi_{\al,\bet}(E)|&= p^{-2ek}\sum_{s,t\in R}\sum_{x,y,z\in E} \chi(s(x\cdot y-\al))\chi(t(x\cdot z - \bet)) \\
                         &= I+II+III,
   \end{align*}
   where $I$ has $s=t=0$, $II$ has $s$ or $t$ equal to zero but not both, and $III$ has $s\neq 0$ and $t\neq 0$.

   For $I$, we see 
   \begin{equation}\label{eqCh3PairsI}
      I=p^{-2ek}\sum_{s=t =0}\sum_{x,y,z \in E}\chi(-0\al)\chi(0\bet)\chi(x\cdot (0y-0z))=p^{-2ek}|E|^3.
   \end{equation}


   With $II$ and $III$, we will use Theorem \ref{singleDot} from Section \ref{secSingDotProd}, which requires that 
   $$ |E| \geq p^{ek}\sum_{i=0}^{e-1} \sqrt{p^{dek+dik-ik}\left(2p^{-ik-k}+1+p^{-ek}\right)}$$
   which is given. 

   For $II$, as either $s=0$ or $t=0$ but not both, we may split the sum as two sums where either $s=0$ or $t=0$,

   \begin{equation}\label{eqCh3PairsII}
      II=p^{-ek}\left(p^{-ek}\sum_{s \in R, t=0}\sum_{x,y,z \in E}\chi(s(x\cdot y-\al))+p^{-ek}\sum_{t \in R, s=0}\sum_{x,y,z \in E}\chi(t(x\cdot z-\bet))\right).
   \end{equation}
   Notice that this is the sum of two $\nu(t)$ from Theorem \ref{singleDot}. Thus by Theorem \ref{singleDot}, we have $II<4|E|^2/p^{2ek}$. 

   Considering $III$, we have as $s,t \neq 0$,  
   \begin{equation}\label{ch3IIIdef}
      |III| \leq \left|p^{-2ek}\sum_{s,t\in R_{e,k}\setminus \{0\}}\sum_{x,y,z\in E}\chi(s(x\cdot y-\al))\chi(t(\bet-x\cdot z))\right|
   \end{equation}
   \begin{align*}
      \intertext{which we then prepare for Cauchy-Schwarz. So,}
            |III|&\leq \left|p^{-2ek}\sum_{x\in E}\sum_{y\in E}\sum_{s\neq 0}\chi(s(x\cdot y-\al))\sum_{t\neq 0}\sum_{z\in E}\chi(t(\bet-x\cdot z))\right| \\
            & \leq p^{-2ek}\sum_{x\in R^d}\left|\sum_{y\in E}\sum_{s\neq 0}\chi(s(x\cdot y-\al))\right|\left|\sum_{t\neq 0}\sum_{z\in E}\chi(t(\bet-x\cdot z))\right|. \\
            \intertext{By Cauchy-Schwarz, }
   \end{align*}
   \begin{equation}\label{ch3PdotProdIIIasSE}
            |III| \leq p^{-2ek}\left(\sum_{x\in R^d}\left|\sum_{s\neq 0}\sum_{y\in E}\chi(s(x\cdot y - \al))\right|^2\right)^{1/2}\left(\sum_{x\in R^d}\left|\sum_{t\neq 0}\sum_{z\in E}\chi(t(x\cdot z - \bet))\right|^2\right)^{1/2}.
   \end{equation}
   Which we relabel as
   \begin{align*}
            |III| & \leq p^{-2ek}III_\al \cdot III_\bet.
   \end{align*}

   Notice that $III_\al$ and $III_\bet$ are similar, so we will examine $III_\al$.
   With $III_\al$ we prepare to apply Cauchy-Schwarz a second time. Note,

   \begin{align*}
      III_\al^2 &= \sum_{x\in R^d}\left|\sum_{s\neq 0}\sum_{y\in E}\chi(s(x\cdot y - \al))\right|^2.\\
   \end{align*}
   So,
   \begin{equation}\label{doubIIIalBound}
                III_\al^2 = \sum_{x\in R^d}\left|\sum_{i=0}^{e-1}\left[ (1)\left(\sum_{s\in [p^i]}\sum_{y\in E}\chi(s(x\cdot y-\al))\right)\right]\right|^2. 
   \end{equation}
   \begin{align*}
                \intertext{Which by Cauchy-Schwarz is}
                III_\al^2 & \leq \sum_{x\in R^d}\sum_{i=0}^{e-1}(1)^2\sum_{i=0}^{e-1}\left|\sum_{s\in [p^i]}\sum_{y\in E}\chi(s(x\cdot y-\al))\right|^2 \\
                &\leq e\sum_{x\in R^d}\sum_{i=0}^{e-1}\sum_{s,s'\in [p^i]}\sum_{y,y'\in E}\chi(s(x\cdot y-\al))\ovl{\chi(s'(x\cdot y'-\al))}. \\
                \intertext{Recall that by Lemma \ref{pBrakStruc}, $s\in [p^i]$ means $s=p^iu$ for some unit $u$. Likewise $s'=p^iv$. This allows us to relabel $s,s'\in [p^i]$ as $p^iu,p^iv\in [p^i]$. Since $\chi$ is an additive character, we may rearrange the terms in the innermost sum to the below:}
                III_\al^2&\leq e\sum_{x\in R^d}\sum_{i=0}^{e-1}\sum_{y,y'\in E}\sum_{p^iu,p^iv\in [p^i]}\chi(p^i(uy-vy')\cdot x)\chi(p^i\al(v-u)). \\
   \end{align*}


   \begin{align*}
      \intertext{We now introduce a factor of $|E|p^{-2ek}$ into the outer-most sum to prepare $III_\al^2$ for Corollary \ref{dotProdRemNu}:}
                III_\al^2&\leq (|E|p^{-2ek})\inv e\sum_{i=0}^{e-1}|E|p^{-2ek}\sum_{x\in R^d}\sum_{y,y'\in E}\sum_{p^iu,p^iv\in [p^i]}\chi(p^i(uy-vy')\cdot x)\chi(p^i\al(v-u)). \\
                \intertext{Which by Corollary \ref{dotProdRemNu} is}
                 &\leq (|E|p^{-2ek})\inv e\sum_{i=0}^{e-1}\left(2|E|^2p^{dek+dik-2ik-k} + |E|^2p^{dek+dik-ik}\right. \left.  + |E|^2p^{dek+dik-ek-ik}\right). \\
                 \intertext{This then simplifies as follows,} 
                III_\al^2 & \leq |E|^{-1}p^{2ek}e\sum_{i=0}^{e-1}|E|^2p^{dek+dik-ik}\left(2p^{-ik-k}  + 1 + p^{-ek}\right) \\
                 & = |E|p^{dek+2ek}e\sum_{i=0}^{e-1}\left(2p^{-k}p^{dik-2ik} + p^{dik-ik} + p^{-ek}p^{dik-ik}\right). \\
    \end{align*}
    Summing each term in $i$ gives,
    \begin{equation}\label{IIIalsqrdleqEq}
                 III_\al^2 \leq |E|p^{dek+2ek}e\left(2p^{-k}\frac{p^{dek-2ek}-1}{p^{dk-2k}-1} + (1+  p^{-ek})\frac{p^{dek-ek}-1}{p^{dk-k}-1}\right). 
   \end{equation}
   Let 
   \begin{align*}
      F &=|E|p^{dek+2ek}e\left(2p^{-k}\frac{p^{dek-2ek}-1}{p^{dk-2k}-1} + (1+  p^{-ek})\frac{p^{dek-ek}-1}{p^{dk-k}-1}\right).\\
   \end{align*}
      By equation \ref{IIIalsqrdleqEq}, we have $III_\al^2 \leq F$.
   Likewise $III_\bet^2 \leq F$. Ergo, 
      $III_\al III_\bet \leq F$.

      Thus as $|III|\leq p^{-2ek}III_\al III_\bet$,
   \begin{equation}\begin{split}\label{ch3PdotProdIII}
                 |III| &\leq  p^{-2ek}F.
   \end{split}\end{equation}

   Putting $I,\ II,\  III$ together, we see that $$I+II+III\leq |E|^3/p^{2ek}+4|E|^2/p^{2ek}+p^{-2ek}F$$
   As we want $I+II+III=|E^3|/p^{2ek}(1+o(1))$, we need show that $II+III=(|E|^3/p^{2ek})o(1)$ when $|E|$ is of sufficient size. Note that by definition of $F$, we have $e|E|p^{dek}/p^{2ek} \leq F/p^{2ek}$. As $4|E|^2/p^{2ek}< e|E|p^{dek}/p^{2ek}$ for any size of $E$ (recall $E\subset R^{d}$ and $|R|=p^{ek}$), we need only have $2p^{-2ek}F \leq |E|^3p^{-2ek}$.

   Let $G=F/|E|$. 
   This gives rise to the inequality,
   \begin{align*}
      |E|^3p^{-2ek} &\geq 2p^{-2ek}F. \\
      \intertext{Which simplifies to}
      |E|^2 & \geq 2G \\
      |E| & \geq \sqrt{2G}.
   \end{align*}
   Thus, as we used Theorem \ref{singleDot} earlier in the proof and it requires 
   $$|E|\geq p^{ek}\sum_{i=0}^{e-1} \sqrt{p^{dek+dik-ik}\left(2p^{-ik-k}+1+p^{-ek}\right)},$$
   we must have
   \begin{align*}
      |E| & \geq \max\left(p^{ek}\sum_{i=0}^{e-1} \sqrt{p^{dek+dik-ik}\left(2p^{-ik-k}+1+p^{-ek}\right)}, \sqrt{2G}\right).
   \end{align*}
\end{proof}


   We now give a few lemmas that will be used later that come from the above proof.
   \begin{defn}
      Let  $$S_{E,\al} := \sum_{s\in R\setminus \{0\}}\sum_{y\in E}\chi(s(x\cdot y - \al)).$$ 
   \end{defn}
   \begin{lemma}\label{pairsSE}
	   Let $d\geq 3, e\geq 5, k\geq 1$. Let $E\subset R_{e,k}^d$.
      Let $$|E| \geq \max\left(p^{ek}\sum_{i=0}^{e-1} \sqrt{p^{dek+dik-ik}\left(2p^{-ik-k}+1+p^{-ek}\right)}, \sqrt{2G}\right).$$ 
      Then,
      \begin{align*}
         & \left(\sum_{x\in E}\left|S_{E,\al}(x)\right|^2\right)^{1/2}\left(\sum_{z\in R_{e,k}}\left|S_{E,\bet}(z)\right|^2\right)^{1/2} \\
         & \qquad \leq |E|p^{dek+2ek}e\left(2p^{-k} \frac{p^{dek-2ek}-1}{p^{dk-2k}-1} +\left(1+p^{-ek} \right)\frac{p^{dek-ek}-1}{p^{dk-k}-1} \right).
      \end{align*}
   \end{lemma}
   \begin{proof}
       From Equation \ref{ch3PdotProdIIIasSE} we get $p^{2ek}|III|$ is less than our left hand side and from Equation \ref{ch3PdotProdIII} we get $p^{2ek}|III|$ is less than our right hand side. As we only relax our bound on $|III|$ between the two equations, the lemma follows.
   \end{proof}
   
   \begin{lemma}\label{pairsSEClean}
      Let $d\geq 3, e\geq 5, k\geq 1$ and $E\subset R_{e,k}^d$.
      Let $$|E| \geq \max\left(p^{ek}\sum_{i=0}^{e-1} \sqrt{p^{dek+dik-ik}\left(2p^{-ik-k}+1+p^{-ek}\right)}, \sqrt{2G}\right).$$ 
      Then,
      $$\left(\sum_{x\in E}|S_{E,\al}(x)|^2\right)^{1/2}\left(\sum_{z\in R_{e,k}}|S_{E,\bet}(z)|^2\right)^{1/2} \leq 7e|E|p^{2dek+ek-dk+k}.$$
   \end{lemma}
   \begin{proof}
      From Lemma \ref{pairsSE}, we have:
      \begin{align*}
         & \left(\sum_{x\in E}\left|S_{E,\al}(x)\right|^2\right)^{1/2}\left(\sum_{z\in R_{e,k}}\left|S_{E,\bet}(z)\right|^2\right)^{1/2} \\
         & \qquad \leq |E|p^{dek+2ek}e\left(2p^{-k} \frac{p^{dek-2ek}-1}{p^{dk-2k}-1} +\left(1+p^{-ek} \right)\frac{p^{dek-ek}-1}{p^{dk-k}-1} \right)
      \end{align*}
      which relaxes to
      \begin{align*}
         & \left(\sum_{x\in E}\left|S_{E,\al}(x)\right|^2\right)^{1/2}\left(\sum_{z\in R_{e,k}}\left|S_{E,\bet}(z)\right|^2\right)^{1/2} 
          \leq 7e|E|p^{2dek+ek-dk+k}
      \end{align*}
   \end{proof}


\subsection{Tree Configurations of Dot Products Proof}


We will now prove Theorem \ref{trees}. We will make use of the following lemma from \cite{steele}.
\begin{lemma}\label{doubSum}
   Let $m$ and $n$ be positive integers. Then for each double sequence $\{c_{jk}:1\leq j\leq m,\ 1\leq k\leq n\}$ and pair of sequence $\{z_j: 1\leq j\leq m\}$ and $\{y_k:1\leq k\leq n\}$ of complex numbers, we have the bound $$\left|\sum_{j=1}^m\sum_{k=1}^nz_jy_kc_{jk}\right|\leq \sqrt{RC}\left(\sum_{j=1}^m|z_j|^2\right)^{1/2}\left(\sum_{k=1}^n|y_k|^2\right)^{1/2}$$ where $R=\max_{1\leq j \leq m} \sum_{k=1}^m c_{jk}$ and $C=\max_{1\leq k \leq m} \sum_{j=1}^m c_{jk}$, the row and column sum maxima.
\end{lemma}

\begin{dff}\label{antenaDff}
   By antenna, we will mean a leaf (a vertex with only one edge) and its corresponding edge. 
\end{dff}
We now prove Theorem \ref{trees}.
\begin{proof}

Let $T$ be as stated in the theorem.
Recall Definition \ref{PiTreeDef}. By Orthogonality we see that
 $$\Pi_{\al}(T) = p^{-nek}\sum_{\substack{x_v \in E \\ v\in V}}\sum_{\substack{s_\eps \in R \\ \eps \in \edgeSetVar}}\prod_{\eps\in \edgeSetVar}\chi(s_\eps(x_{\eps_1}\cdot x_{\eps_2} - \al_{\eps})).$$

Let $R_\ell$ be the sum of the sums with $\ell$ many $s_j\neq0$.
In the case that all $s_j=0$, we get that 

\begin{align*}
   |R_0| &= p^{-nek}\sum_{\substack{x_v \in E \\ v\in V}}\sum_{\substack{s_\eps = 0 \\ \eps \in \edgeSetVar}}\prod_{\eps\in \edgeSetVar}\chi(s_\eps(x_{\eps_1}\cdot x_{\eps_2} - \al_{\eps})) \\
         &= p^{-nek}\sum_{\substack{x_v \in E \\ v\in V}}\sum_{\substack{s_\eps = 0 \\ \eps \in \edgeSetVar}}\prod_{\eps\in \edgeSetVar}1 \\
         &= p^{-nek}|E|^{|V|} \\
         &= \frac{|E|^m}{p^{nek}}.
\end{align*}

For $R_1$, we first consider an arbitrary configuration $R_{\vartreeNumII}$ with only one $s_j\neq 0$.  Let edge $\varedgeII = (v_1,v_2)$ correspond with this choice of $s_j$. We may factor the $\chi(s_j(v_1\cdot v_2 - \al_{\varedgeII}))$ term out as any other character with $v_1$ or $v_2$ in it must be equal to 1. This gives, 
\begin{align*}
   R_\vartreeNumII & \leq p^{-nek}\sum_{v_2\in E}\left(\sum_{v_1\in E}\sum_{s_j\neq 0}\chi(s_j(v_1\cdot v_2 - \al_{\varedgeII}))\right)\sum_{\substack{x_v \in E \\ v\in V\setminus \{v_2,v_1\}}}\sum_{\substack{s_\eps = 0 \\ \eps \in \edgeSetVar\setminus \{\varedgeII\}}}\prod_{\eps\in \edgeSetVar\setminus \{\varedgeII\}}\chi(0(x_{\eps_1}\cdot x_{\eps_2} - \al_{\eps})) \\
          & \leq p^{-nek}\sum_{v_2\in E}\left(\sum_{v_1\in E}\sum_{s_j\neq 0}\chi(s_j(v_1\cdot v_2 - \al_{\varedgeII}))\right)|E|^{|V|-2}. \\
           \intertext{ By Corollary \ref{dotProdRemSimp} } 
   R_\vartreeNumII & \leq p^{-nek+ek}\sqrt{6+3e}|E|p^{dek-dk/2+ek/2+k/2}|E|^{|V|-2}. \\
\end{align*}
As we have $n$  choices for our one edge $\varedgeII = (v_1,v_2)$ with $s_j\neq 0$,
this gives $$|R_1| \leq p^{-nek}\sqrt{6+3e}n|E|^{m-1}p^{dek-dk/2+3ek/2+k/2}.$$

For the case of $R_\ell$ with $\ell>1$, we consider first an arbitrary configuration $R_{\vartreeNumIII}$ with $\ell$ many $s_j\neq 0$. Let $\edgeSetVar_{\neq}$ be the set of edges with $s_j\neq 0$ and $\edgeSetVar_=$ the set with $s_j=0$. Then as we have a tree and $\ell>0$, there must exist two antennae with edges in $\edgeSetVar_{\neq}$, let there vertices be $a_1,a_2$ and edges be $e_1=(a_1,b_1)$ and $e_2=(a_2,b_2)$ with $a_1\cdot b_1 = \al_{ab_1}$ and $a_2\cdot b_2 = \al_{ab_2}$. Then 
\begin{align*}
   |R_\vartreeNumIII| &= p^{-nek}\sum_{b_1\in E}\left(\sum_{s_{ab_1}\in R\setminus \{0\}} \sum_{a_1\in E}\chi_{s_{ab_1}}(a_1\cdot b_1 - \al_{ab_1})\right)\\
                     & \qquad \tm \sum_{b_2\in E}\left(\sum_{s_{ab_2}\in R\setminus \{0\}} \sum_{a_2\in E}\chi_{s_{ab_w}}(a_2\cdot b_2 - \al_{ab_2})\right)F(b_1,b_2) \\
\end{align*}
where $F$ represents the rest of the sum $$F(b_1,b_2) =  \sum_{\substack{x_v \in E \\ v\in V\setminus \{b_2,b_1\}}}\sum_{\substack{s_\eps = 0 \\ s_\eps' \neq 0 \\ \eps \in \edgeSetVar_= \\ \eps' \in \edgeSetVar_{\neq}\setminus \{e_1,e_2\}}}\prod_{\eps\in \edgeSetVar\setminus \{e_1,e_2\}}\chi(s_\eps(x_{\eps_1}\cdot x_{\eps_2} - \al_{\eps})).$$

For Lemma \ref{doubSum}, let $j=b_1$, $k=b_2$, and $z_j=\sum_{s_{ab_1}\in R\setminus \{0\}} \sum_{a_1\in E}\chi_{s_{ab_1}}(a_1\cdot b_1 - \al_{ab_1})$, and $y_k = \sum_{s_{ab_2}\in R\setminus \{0\}} \sum_{a_2\in E}\chi_{s_{ab_w}}(a_2\cdot b_2 - \al_{ab_2})$, and $c_{jk}=F$
(viewing $F$ as a matrix by $F_{b_1b_2}=F(b_1,b_2)$).
By Lemma \ref{pairsSEClean} we get an upper bound of $7e|E|p^{2dek-dk+ek+k}$ for our factored out edges. By Lemma \ref{doubSum} we have a bound on $F$ by $\sqrt{RC}$ where $R$ and $C$ are the row and column sum maxima of $F(b_1,b_2)$. That is $R=\max_{b_1\in E}\sum_{b_2\in E}|F(b_1,b_2)|$. Together, we have

\begin{align*}
   |R_\vartreeNumIII| & \leq p^{-nek}7e|E|p^{2dek-dk+ek+k}\cdot \sqrt{RC}.
\end{align*}
 By the triangle inequality, and that characters are bounded above by $1$, may bound the product inside of $F$ by $1$, giving 
\begin{align*}
   |F(b_1,b_2)| &\leq \sum_{\substack{x_v \in E \\ v\in V\setminus \{b_2,b_1,a_1,a_2\}}}\sum_{\substack{s_\eps = 0 \\ s_\eps' \neq 0 \\ \eps \in \edgeSetVar_= \\ \eps' \in \edgeSetVar_{\neq}\setminus \{e_1,e_2\}}} 1. \\
\end{align*}

Then each vertex that $F$ sums over contributes a factor of $|E|$, each edge with $s_j=0$ that $F$ sums over contributes a factor of $1$, and each edge with $s_j\neq 0$ contributes a factor of $p^{ek}$, giving that $F\leq |E|^{m-4}p^{ek(\ell-2)}$ and so $R \leq |E|^{m-3}p^{ek(\ell-2)}$. Likewise, $C\leq |E|^{m-3}p^{ek(\ell-2)}$. This gives
$$|R_\vartreeNumIII| \leq p^{-nek}7e|E|^{m-2}p^{ek(\ell-2)}p^{2dek-dk+ek+k}$$ and so since we have ${ n \choose \ell}$ ways to choose our $R_\vartreeNumIII$, $$R_\ell \leq p^{-nek}{n \choose \ell}7e|E|^{m-2}p^{ek(\ell-2)}p^{2dek-dk+ek+k}.$$

As we want to find the size of $E$ such that $R_0\geq \sum_{i=1}^n R_n$, we have
\begin{align*}
   \frac{|E|^m}{p^{nek}} & \geq p^{-nek}n\sqrt{6+3e}|E|^{m-1}p^{dek-dk/2+3ek/2+k/2}+\sum_{\ell=2}^n p^{-nek}{n \choose \ell}7e|E|^{m-2}p^{ek(\ell-2)}p^{2dek-dk+ek+k} \\
   |E|^2 & \geq n\sqrt{6+3e}|E|p^{dek-dk/2+3ek/2+k/2}+7ep^{-2ek}(p^{ek}+1)^np^{2dek-dk+ek+k} \\
   |E| & \geq \sqrt{n\sqrt{6+3e}|E|p^{dek-dk/2+3ek/2+k/2}+7ep^{-2ek}(p^{ek}+1)^np^{2dek-dk+ek+k}} \\
   |E| & \geq \sqrt{2}\max\left\{\sqrt{n\sqrt{6+3e}|E|p^{dek-dk/2+3ek/2+k/2}},\ \sqrt{7ep^{-2ek}(p^{ek}+1)^np^{2dek-dk+ek+k}}\right\}. \\
\end{align*}
For the first branch of the $\max$,
\begin{align*}
   |E| & \geq \sqrt{2n\sqrt{6+3e}|E|p^{dek-dk/2+3ek/2+k/2}} \\
   |E| & \geq 2n\sqrt{6+3e}p^{dek-dk/2+3ek/2+k/2}. \\
\end{align*}
For the second,
\begin{align*}
   |E| & \geq \sqrt{14}ep^{-ek}(p^{ek}+1)^{n/2}p^{dek-dk/2+ek/2+k/2} \\
       & \geq \sqrt{14}ep^{dek-dk/2+nek/2+n/2-ek/2+k/2}.
\end{align*}
When $n\geq 4$, $|E| \geq \sqrt{14}enp^{dek-dk/2+nek/2+n/2-ek/2+3k/2}.$

\end{proof}

\subsection{Corollaries of Tree Configurations Result}
\subsubsection{K-Chains}

\begin{defn}
   A $k$-chain is a $(k+1)$-tuple of distinct points $(x_1,\ldots, x_{k+1})\in E^{k+1}$ such that $x_j\cdot x_{j+1}=\al_j$ for every $1\leq j \leq k$. 
\end{defn}

\begin{cor}
   Let $d\geq 3, e\geq 5, k\geq 1, n\geq 4$. Let $E \subset R_{e,k}^d$. Let $\al=(\al_1,\ldots, \al_n)$.  Then the number of $k$-chain configurations, $\Pi_\al(K)$, has $\Pi_\al(K)\leq \frac{2|E|^{k+1}}{p^{nek}} $ whenever $|E|>\sqrt{14}enp^{dek-dk/2+nek/2+n/2-ek/2+3k/2}$.
\end{cor}
   \begin{proof}
As this relation $x_j\cdot x_{j+1}=\al_j$ forms no cycles, a $k$-chain is a type of forest. Note we have $n+1$ many elements and $n$ many relations. Let $K$ be the graph formed by the $n+1$ elements as vertices and the relations $x_j\cdot x_{j+1}=\al_j$ for the edges. Then by Theorem \ref{trees}, we have $$\Pi_\al(K) \leq \frac{2|E|^{n+1}}{p^{nek}} $$ whenever $|E| > \sqrt{14}enp^{dek-dk/2+nek/2+n/2-ek/2+3k/2}$.
   \end{proof}

\subsubsection{Stars}

\begin{defn}
   A star is a $(k+1)$-tuple of distinct points $(c,x_1,\ldots, x_{k})\in E^{k+1}$ such that $c\cdot x_j=\al_j$ for every $1\leq j \leq k$ for some given $\{\al_j\}_{1\leq j\leq k}$. 
\end{defn}

\begin{cor}
   Let $d\geq 3, e\geq 5, k\geq 1, n\geq 4$. Let $E \subset R_{e,k}^d$. Let $\al=(\al_1,\ldots, \al_n)$.  Then the number of star configurations $\Pi_al(K)$ has $\Pi_\al(K)\leq \frac{2|E|^{k+1}}{p^{nek}} $ whenever $|E|>\sqrt{14}enp^{dek-dk/2+nek/2+n/2-ek/2+3k/2}$.
\end{cor}
\begin{proof}
 As this relation $c\cdot x_j=\al_j$ forms no cycles, a star is a type of forest. Note we have $n+1$ many elements and $n$ many relations. Let $K$ be the graph formed by the $n+1$ elements as vertices and the relations $c\cdot x_j=\al_j$ for the edges. Then by Theorem \ref{trees}, we have $$\Pi_\al(K) \leq \frac{2|E|^{n+1}}{p^{nek}} $$ whenever $|E| > \sqrt{14}enp^{dek-dk/2+nek/2+n/2-ek/2+3k/2}$.
\end{proof}

\subsubsection{Trees in $\ZZ_{p^l}$}
Let $R=\ZZ_{p^e}$, then $R=R_{e,1}$ and so our $k=1$. Let $q=p^e$. This gives that for a dot product configuration tree $T$ with $m$ many elements and $n$ many relations, $$\Pi(T) \leq \frac{2|E|^m}{p^{nek}}$$ whenever 
$|E| \geq \sqrt{14}enp^{de-d/2+ne/2+n/2-e/2+3/2}$ or in terms of $q$, $|E|\geq \sqrt{14}enq^{\frac{2d-d/e}{2}+n/2-1/2}p^{n/2+3/2}.$

\subsubsection{Trees in $\FF_{p^l}$}

Let $R=\FF_{p^k}$, then $R=R_{1,k}$ and so our $e=1$. Let $q=p^k$. This gives that for a dot product configuration tree $T$ with $m$ many elements and $n$ many relations, $$\Pi(T) \leq \frac{2|E|^m}{p^{nek}}$$ whenever $|E| \geq \sqrt{14}np^{dk-dk/2+nk/2+n/2}$ or in terms of $q$, $|E|\geq \sqrt{2}nq^{d/2+n/2}p^{n/2}.$

\subsection{Inverse Row-Matrix Problem}

\begin{proof}
    Let $d\geq 3,\ e\geq 5,\ k\geq 1,\ n\geq 4$. Let $R=R_{e,k}$. Let $E\subset R^d$.
As a star is a type of tree, we know Theorem \ref{trees} holds for stars. 
Also note that a star is $\{(c,x_1,\ldots, x_{n}),\ c,x_j \in E \mid c\cdot x_j = \al_j\}$ which has $n+1$ many vertices and $n$ many edges. As a system of linear equations, this is
\[
   \left\{
\begin{array}{ccccccccc}
   c_1x_{11} &+& c_2x_{12} &+& \cdots &+& c_{d}x_{1d} &=& \al_1 \\
   c_1x_{21} &+& c_2x_{22} &+& \cdots &+& c_{d}x_{2d} &=& \al_2 \\
   \vdots    &+& \vdots    &+& \ddots &+& \vdots      &=& \vdots \\
   c_1x_{n1} &+& c_2x_{n2} &+& \cdots &+& c_dx_{nd}   &=& \al_d
\end{array}
\right.
\]

This is also the matrix equation

$$c\mat{x_1 & \cdots &x_{n}} = \al. $$ 

Relabeling, we have $$aA=b$$ where $a\in E$ and the columns of $A$ are in $E$ and $b$ is given. Let $K$ be as from Stars. By Theorem \ref{trees}, whenever $|E|> \sqrt{14}enp^{dek-dk/2+nek/2+n/2-ek/2+3k/2}$, then $\left|\Pi_\al(K)\right| \leq \frac{2|E|^{n+1}}{p^{nek}}$. 

By Definition \ref{PiTreeDef},
\begin{align*}
   \Pi_\al(K) &= \{ x \in E^{n+1} : \forall \al_{ij} \in \al, x_i\cdot x_j = \al_{ij}\} \\
              &= \{(c,x_1,\ldots x_n) \in E^{n+1} : \forall \al_i \in \al, c\cdot x_i = \al_i\} \\
              &\cong \{ (c,[x_1 \ldots x_n]) \in E\tm E^n : c[x_1 \ldots x_n] = \al\}. \\
              \intertext{where $\cong$ means isomorphic as sets. By our relabeling, this is,}
              \Pi_\al(K) &\cong \{ (a,A) : a\in E, A\in E^n, aA=b\}.
\end{align*}

So, the number of solutions $$S=\left|\left\{(a,A)\mid a\in E, A\in E^{n},\  aA=b\right\}\right|$$ has 
$$ S \leq \frac{2|E|^{n+1}}{p^{nek}} $$
whenever $|E|> \sqrt{14}enp^{dek-dk/2+nek/2+n/2-ek/2+3k/2}$.
\end{proof}
\bibliographystyle{amsplain}


\end{document}